\documentclass[pdflatex,sn-mathphys]{sn-jnl}

\usepackage{bm}
\usepackage[group-separator={,}]{siunitx}
\usepackage[T1]{fontenc} 

\jyear{2022}%

\theoremstyle{thmstyleone}%
\theoremstyle{thmstyletwo}%

\theoremstyle{thmstylethree}%

\begin{document}

\title[DD-CPM library for PDEs on surfaces]{A closest point method library for PDEs on surfaces with parallel domain decomposition solvers and preconditioners}
%

\author*[1]{\fnm{Ian C.T.} \sur{May}}\email{imay1@ucsc.edu}
\author[2]{\fnm{Ronald D.} \sur{Haynes}}\email{rhaynes@mun.ca}
\author[3]{\fnm{Steven J.} \sur{Ruuth}}\email{sruuth@sfu.ca}

\affil*[1]{\orgdiv{Department of Applied Mathematics}, \orgname{University of California Santa Cruz}, \orgaddress{\city{Santa Cruz}, \postcode{95064}, \state{California}, \country{USA}}}

\affil[2]{\orgdiv{Department of Mathematics}, \orgname{Memorial University of Newfoundland}, \orgaddress{\city{St. John's}, \postcode{A1C 5S7}, \state{Newfoundland}, \country{Canada}}}

\affil[3]{\orgdiv{Department of Mathematics}, \orgname{Simon Fraser University}, \orgaddress{\city{Burnaby}, \postcode{V5A 1S6}, \state{British Columbia}, \country{Canada}}}

\abstract{The DD-CPM software library provides a set of tools for the discretization and solution of problems arising from the closest point method (CPM) for partial differential equations on surfaces. The solvers are built on top of the well-known PETSc framework, and are supplemented by custom domain decomposition (DD) preconditioners specific to the CPM. These solvers are fully compatible with distributed memory parallelism through MPI. This library is particularly well suited to the solution of elliptic and parabolic equations, including many reaction-diffusion equations. The software is detailed herein, and a number of sample problems and benchmarks are demonstrated. Finally, the parallel scalability is measured.}


\maketitle
\section{Introduction}
\label{sec:intro}
The numerical solution of PDEs intrinsic to surfaces presents interesting challenges over their flat-space analogs. The closest point method (CPM) is an approach to discretize general surface-intrinsic PDEs over a wide class of surfaces. The software library presented here provides an MPI parallel implementation of the CPM built on top of PETSc \cite{petsc-efficient}. This software is particularly well suited to solving surface intrinsic reaction-diffusion systems, generically of the form
\begin{equation}
  \left(\frac{\partial}{\partial t} - \mu_i\Delta_{\mathcal{S}}\right)u_i = f_i(t,\bm{x},\bm{u}),
  \label{eq:rdsys}
\end{equation}
where $\bm{u}$ are the different species in the system, $\mu_i$ are the diffusion coefficients, and $f_i$ are the nonlinear reaction terms coupling the systems together.

This software provides access to the extensive suite of tools within PETSc for use with the CPM, and additionally defines custom Schwarz-type domain decomposition solvers/preconditioners for equations of the form
\begin{equation}
  \left(c-\Delta_{\mathcal{S}}\right)u = f,
  \label{eq:shiftPoisson}
\end{equation}
where $c\in\mathbb{R}^+$ is a positive constant, and $\Delta_{\mathcal{S}}$ is the Laplace-Beltrami operator intrinsic to an embeddable surface $\mathcal{S}\subset\mathbb{R}^d$. These custom domain decomposition (DD) methods were first presented in \cite{May:DD25}, and subsequently expanded upon in \cite{May:SISc}. Solving equation~(\ref{eq:rdsys}) semi-implicitly requires the solution of many elliptic systems. The provided custom DD preconditioners are well suited to these linear systems. Additionally, the flexibility of PETSc is maintained, and the use of other preconditioners and solvers within PETSc is straightforward.

To date, most software for the CPM has been limited to the use of direct solvers on shared memory machines. The DD-CPM software detailed here provides the first known implementation of this method that can leverage larger distributed memory machines. This facilitates the application of the CPM to larger and more complicated problems. For elliptic PDEs, the dependence on direct solvers that has been present in most CPM implementations has also been lifted through the introduction of custom DD preconditioned Krylov methods compatible with the CPM. For reaction-diffusion equations, this software allows the use of fine grids to reliably capture either complex surface geometry and/or sharp features in the generated solutions. Indeed, the software has been kept general to allow many other equations to be posed and solved.

\subsection{Existing software}
\label{sec:intro:exist}
There is little existing software for the CPM. An important resource is the {\tt cp\_matrices} repository \cite{cpmatrices}, hosted by Prof. Colin Macdonald (coauthor of references \cite{CBM:ICPM,CBM:Eig,CBM:RDonPC,Maerz/Macdonald:cpfunctions,CBM:LSE,Chen:MG}). This repository consists of MatLab and Python codes, and is mostly restricted to serial execution. For problems outside our framework, we suggest to start with the {\tt cp\_matrices} repository.

A large motivator for developing CPM specific DD methods was the experience that algebraic preconditioners did not perform well in all cases, and required a great deal of problem specific tuning. However, there may be cases where the built-in PETSc preconditioners perform well with little cost per iteration \cite{MK:AMG}. The back-end preconditioners and solvers used with the DD-CPM library can be substituted with PETSc options easily, and without re-compiling the library. An important goal of the development of this software is to maintain a high degree of flexibility and interoperability with other software packages.

\section{Review of core methods}
\label{sec:rev}
Before discussing the software, we will briefly review the closest point method (CPM), and the (optimized) restricted additive Schwarz (ORAS) domain decomposition method. To keep these reviews concrete and simple they will be constrained to equation \eqref{eq:shiftPoisson}; the consideration of more complicated equations will be deferred to Section \ref{sec:running:user}.

\subsection{Closest point method}
\label{sec:rev:cpm}
Methods for the numerical solution of surface intrinsic PDEs generally take one of two approaches: discretize the surface itself, or find a solution in a higher dimensional embedding space. Surface parametrization methods \cite{FloaterHormann:Para} are very efficient and allow the use of familiar discretizations, but are limited to simple surfaces where the parametrization is known.  With these methods the user may need to contend with coordinate singularities. Finite element methods acting on a triangulation of the surface \cite{DziukElliot} yield sparse symmetric systems for model equation \eqref{eq:shiftPoisson}, but are sensitive to the quality of the triangulation. Level set methods \cite{Cheng:LSM} embed the surface in a higher dimensional flat space, treating the surface implicitly, however the requisite artificial boundary conditions on the computational domain and the treatment of open surfaces are non-trivial.

The CPM \cite{SJR:CPM,CBM:ICPM} is an embedding method that uses an implicit representation of the surface similar to the level set method. When implemented, the CPM is posed only over a small region of the embedding space, $\mathbb{R}^d$, near the surface. This dramatically reduces the number of unknowns in the discretization. The CPM has been applied to triangulated surfaces \cite{CBM:LSE}, surfaces of mixed codimension \cite{SJR:CPM}, moving surfaces \cite{Argy:MoveSurf}, and even point cloud domains \cite{CBM:RDonPC}. The DD-CPM library supports triangulated surfaces directly. Surfaces of mixed codimension and point cloud domains are both supported by allowing user defined closest point functions to be used. This software does not currently support moving surfaces.

In this section, the CPM is described for equation \eqref{eq:shiftPoisson} over a simple closed surface to keep the presentation brief. The CPM utilizes familiar Cartesian discretizations of differential operators for surface PDEs by first extending the solution off the surface to a narrow tube in the embedding space. The extended solution is formed to be constant in the surface normal direction, and is obtained through the closest point function 
\begin{eqnarray}
  cp_{\mathcal{S}}\colon\mathbb{R}^d&\rightarrow&\mathcal{S} \label{eq:cpFunc} \\
  {\bf x}&\mapsto&\underset{{\bf y}\in\mathcal{S}}{arg\,min} \lvert\lvert{\bf x}-{\bf y}\rvert\rvert_2, \nonumber
\end{eqnarray}
which maps each point in the embedding space, ${\bf x}\in\mathbb{R}^d$, to the closest point on the surface. The extension operator $E$ is defined as the composition  of a function $f\colon\mathcal{S}\rightarrow\mathbb{R}$ with $cp_\mathcal{S}$ over the embedding space, i.e.,  $\left(Ef\right)({\bf x}):=f\left(cp_{\mathcal{S}}({\bf x})\right)$ for ${\bf x}\in\mathbb{R}^d$. Importantly, the restriction of the Laplacian of an extended function, $\Delta Ef$, to the surface $\mathcal{S}$, is equivalent to the Laplace-Beltrami operator acting on that surface bound function, $\Delta_{\mathcal{S}}f$ \cite{SJR:CPM,CBM:ICPM,Maerz/Macdonald:cpfunctions}. The core idea of the CPM lies in discretizing $\Delta E$ over an appropriate domain in the embedding space instead of discretizing $\Delta_{\mathcal{S}}$ over the surface directly.

\begin{figure}[htbp]
  \centering
  \begin{minipage}[c]{0.45\linewidth}
    \includegraphics[width=0.95\textwidth]{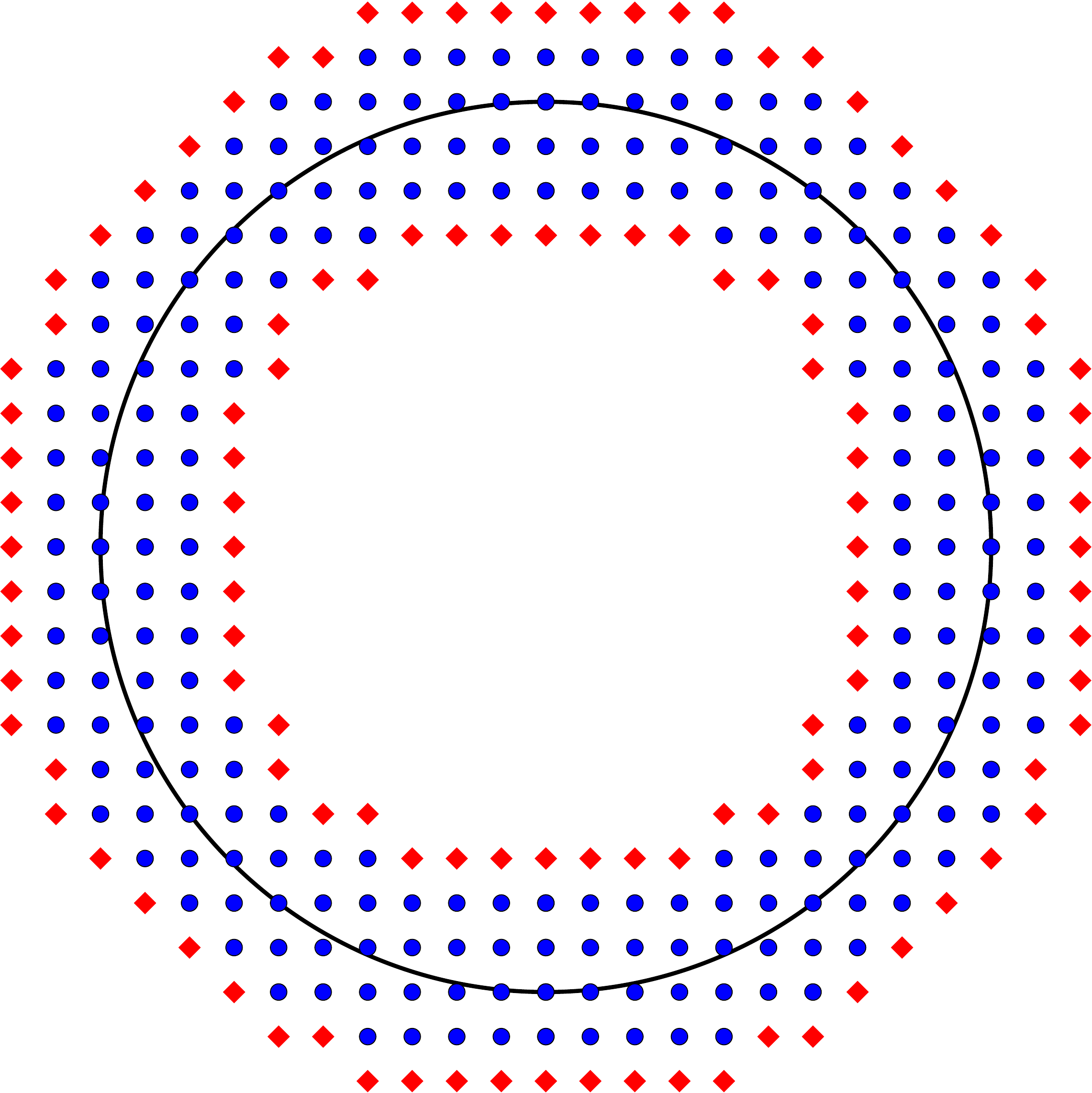}
  \end{minipage}
  \begin{minipage}[c]{0.54\linewidth}
    \includegraphics[width=0.95\textwidth]{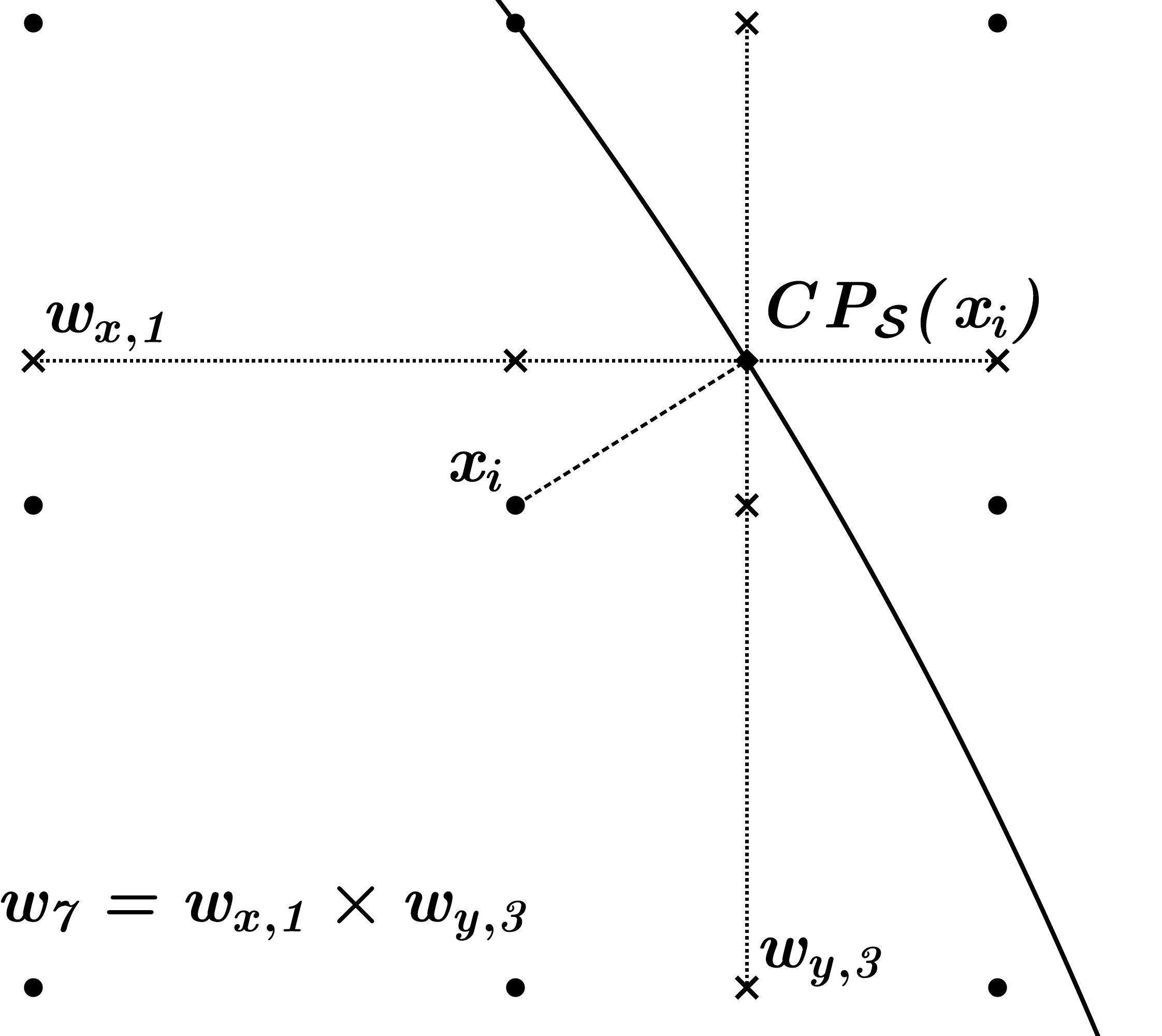}
  \end{minipage}
  \caption{The left panel shows the active nodes $\Sigma_A$ (blue circles) and ghost nodes $\Sigma_G$ (red diamonds) for a circle embedded in $\mathbb{R}^2$. The right panel shows a bi-quadratic extension stencil for the node marked $x_i$. The values at all nine active nodes contribute to the value interpolated to $cp_\mathcal{S}(x_i)$, and thus the extension back to $x_i$.}
  \label{fig:compDom}
\end{figure}

A structured grid with spacing $h$ is placed over the embedding space, and points near the surface are collected into the set of active nodes $\Sigma_A$. The points neighboring the members of $\Sigma_A$ are gathered into the set of ghost nodes $\Sigma_G$. The closest point function \eqref{eq:cpFunc} is continuous on a region of $\mathbb{R}^d$ that lies within a distance $\kappa_\infty^{-1}$ of the surface, where $\kappa_\infty$ is an upper bound on the curvatures of $\mathcal{S}$ \cite{Chu:Vari}. The extent of the active nodes in $\Sigma_A$ will be determined by the extension operator, defined next, and it is assumed that $h$ is chosen small enough that all nodes in $\Sigma_A$ and $\Sigma_G$ lie within a distance $\kappa_\infty^{-1}$ of the surface. The active and ghost nodes for a circle embedded in $\mathbb{R}^2$ can be seen in the left panel of Figure \ref{fig:compDom}.

Given a function $\widetilde{f}$ sampled on $\Sigma_A$, the discrete extension operator ${\rm {\bf E}}$ produces a function on $\Sigma_A$ that is constant in the surface normal direction and coincides with $\widetilde{f}$ on the surface. However, for ${\bf x}_i\in\Sigma_A$ the closest point $cp_\mathcal{S}({\bf x}_i)$ will generally not be another grid point. The function value $\widetilde{f}\left(cp_\mathcal{S}({\bf x}_i)\right)$ is therefore found by interpolating $\widetilde{f}$ over the points in $\Sigma_A$ near $cp_\mathcal{S}({\bf x}_i)$. We use tensor product barycentric Lagrange interpolation \cite{Tref:bary} of degree $p$, and thus $\Sigma_A$ contains the union of all $(p+1)^d$ cubes of grid points needed to perform this extension. The interpolation weights are independent of the data being interpolated, and ${\rm {\bf E}}$ can be written as a matrix where each row holds the interpolation weights for a given node in the computational domain (consisting of $\Sigma_A$ and $\Sigma_G$). An example extension stencil can be seen in the right panel of Figure \ref{fig:compDom}.

The ambient Laplacian, $\Delta$, over the embedding space is discretized by standard second-order centered finite differences, denoted $\Delta^h$. Along the edges of the tube of active nodes in $\Sigma_A$ there will be incomplete finite difference stencils. The ghost nodes in $\Sigma_G$ complete these stencils.

Finally, the discrete form of $\Delta_\mathcal{S}$ will be given as the composition of $\Delta^h$ and ${\rm {\bf E}}$. To apply the approach to implicit time stepping of diffusive problems \cite{CBM:ICPM}, and eigenvalue problems \cite{CBM:Eig}, we use a stabilized discretization given by
\begin{equation}
  \Delta_{\mathcal{S}}^h = -\frac{2d}{h^2}{\rm {\bf I}} + \left(\Delta^h + \frac{2d}{h^2}{\rm {\bf I}}\right){\rm {\bf E}},
  \label{eq:lapBel}
\end{equation}
where the removal of the diagonal elements from $\Delta^h$ avoids redundant self-interpolation.

\subsection{Domain decomposition}
\label{sec:rev:dd}
Iterative methods for the solution of large linear systems are attractive for many reasons. Critical to this work, they may have greatly reduced time and memory requirements over direct solvers, and they are easy to parallelize on distributed memory machines. Domain decomposition (DD) methods seek to replace the solution of one large problem with the repeated solution of several smaller problems, and are particularly well suited to elliptic PDEs \cite{DoleanNataf,ToselliWidlund}. These methods can be used as iterative solvers on their own, or embedded within a Krylov solver as a preconditioner \cite{DoleanNataf}. The DD-CPM library implements Restricted Additive Schwarz (RAS) and Optimized Restricted Additive Schwarz (ORAS) solvers and preconditioners that respect the unique needs of the CPM.

For the CPM it will be beneficial to write the continuous formulation of the (O)RAS methods with respect to the surface intrinsic PDE \eqref{eq:shiftPoisson}. Consider splitting the global domain $\mathcal{S}$ into $N_S$ disjoint subdomains $\widetilde{\mathcal{S}}_j$. Each disjoint subdomain is then grown to form an overlapping set of subdomains $\mathcal{S}_j$. The boundary of each subdomain, $\partial\mathcal{S}_j$, is split into parts lying in the nearby disjoint subdomains denoted as $\Gamma_{jk} \equiv \partial\mathcal{S}_j\cap \widetilde{\mathcal{S}}_k$. Given an initial guess for the global solution $u^{(0)}$, defined at least on the artificial boundaries $\Gamma_{jk}$, we may solve the local problems
\begin{equation}
  \begin{cases}
    \left(c-\Delta_\mathcal{S}\right)u_j = f,&\quad {\rm in}~\mathcal{S}_j,\\
    \mathcal{T}_{jk}u_j = \mathcal{T}_{jk}u^{(n)},&\quad {\rm on}~\Gamma_{jk},~\forall~k.
  \end{cases}
  \label{eq:ddIter}
\end{equation}
After the local solutions $u_j$ are found on each overlapping subdomain, a new global solution may be formed with respect to the disjoint subdomains as
\begin{equation}
  u^{(n+1)} = \sum\limits_j \left.u_j\right\rvert_{\widetilde{\mathcal{S}}_j}.
\end{equation}
Then, with an updated solution on the artificial boundaries, the local problems may be solved again.  The iteration then continues.  The boundary operators $\mathcal{T}_{jk}$ transmit data between the local problems and are thus referred to as \textit{transmission operators}. The DD-CPM software considers transmission operators in one of two forms:
\begin{eqnarray}
  \mathcal{T}_{jk} &=& identity, \\
  \mathcal{T}_{jk} &=& \left(\frac{\partial}{\partial\hat{\rm{\bf q}}_{jk}} + \alpha\right),
  \label{eq:transRob}
\end{eqnarray}
where $\hat{\rm{\bf q}}_{jk}$ is the unit conormal vector along $\Gamma_{jk}$. The first option enforces Dirichlet conditions on the local problems, while the second enforces Robin conditions.   The Robin conditions also provide a parameter $\alpha$ which may be tuned to accelerate convergence. 

ORAS solvers and preconditioners are constructed by discretizing equation \eqref{eq:ddIter} via the CPM. Finding discretizations of the transmission operators $\mathcal{T}_{jk}$, such that the DD iteration converges to the discrete single domain solution, requires care. A full discussion of the requisite considerations, and guidance on setting the parameter $\alpha$, can be found in \cite{May:DD25,May:SISc}.

\section{Overview of the DD-CPM Software}
\label{sec:soft}
The DD-CPM software is hosted as a public git repository on BitBucket. Users may obtain the current stable version from the release branch by calling:
\begin{verbatim}
git clone -b release \
https://mayianm@bitbucket.org/mayianm/DD-CPM.git DD-CPM
\end{verbatim}
Similarly, one can obtain the software as it existed at the time of publication by cloning as above, and then calling:
\begin{verbatim}
git checkout NumAlgoVersion
\end{verbatim}

Detailed installation instructions can be found in the \texttt{README.md} file present in the top level of the repository.

\subsection{Specifying your own equations}
\label{sec:running:user}

Considering the \texttt{DDCPGrayScott} example, the essential ingredients in the driver code can be identified. A \texttt{ProblemDefinition} object, defined and setup on lines $62-71$, holds all of the settings relevant to the setup of all other objects. An object of the \texttt{CPPostProc} class (line $73$) writes out all data files produced by the solver. The \texttt{GridFunc} and \texttt{DiffEq} objects created on lines $77-82$ define the equation being solved. The usage of these classes is discussed in more detail below. The \texttt{CPMeshGlobal} class (used on line $84$) defines the global computational domain and internally manages the partitioning of this domain into subdomains for the DD solvers to use. The \texttt{ProblemGlobal} object on line $86$ builds the global system defined by the \texttt{DiffEq} object over the domain specified by the \texttt{CPMeshGlobal} object, and if DD solvers are desired, will also build the subproblems and manage communication internally. Finally, calling either the \texttt{solveTransient} method (line $91$), or the \texttt{solveStationary}  method (shown in the DDCPPoisson example) on the global problem will solve the defined equations.

\subsubsection{Differential equations}
\label{sec:running:user:diff}
The class \texttt{DiffEq} provides a container that combines one or more \texttt{GridFunc} objects, corresponding to the involved differential operators, with any forcing functions making up the differential equation to be solved. For equations without time dependence, the assumed structure is $\mathcal{L}u = f(x)$ where $\mathcal{L}$ is a differential operator expressed as a \texttt{GridFunc} object (discussed subsequently), and $f(x)$ is a forcing function dependent on the spatial variable. For instance line $67$ in the \texttt{DDCPPoisson} example is:
\begin{lstlisting}[language=C++]
  DiffEq equation(pd,&shiftLap,&forcingFunc);
\end{lstlisting}
which encodes equation \eqref{eq:shiftPoisson} by applying the grid function \texttt{shiftLap} to the left side of the equation, and filling the forcing term from \texttt{forcingFunc} defined at the top of the driver file.

For time dependent equations the assumed structure is $\left(\partial/\partial t + \mathcal{L}\right)u = f(x,t,u)$, where $\mathcal{L}$ is a spatial differential operator represented by a \texttt{GridFunc} object, and $f$ is the forcing function which can now depend on the solution $u$ and the time $t$ in addition to spatial variable $x$. Lines $68$ and $69$ in the \texttt{DDCPHeat} example are:
\begin{lstlisting}[language=C++]
  GridFunc laplacian("LaplacianSecondOrder",
                     pd.mesh.dim,pd.mesh.delta,-1.0);
  DiffEq equation(pd,&laplacian,&icHeat,&rhsHeat);
\end{lstlisting}
where \texttt{laplacian} is a \texttt{GridFunc} object encoding the (negative) Laplace-Beltrami operator. Note now that in addition to the spatial operator and the grid function, the \texttt{DiffEq} object now requires a function specifying the initial condition.

Finally, a \texttt{std::vector} of grid functions can also be supplied to specify multicomponent equations. Consider lines $77-82$ in the \texttt{DDCPGrayScott} example:
\begin{lstlisting}[language=C++]
  GridFunc lapU("LaplacianSecondOrder",
                pd.mesh.dim,pd.mesh.delta,-gsPar_Du);
  GridFunc lapV("LaplacianSecondOrder",
                pd.mesh.dim,pd.mesh.delta,-gsPar_Dv);
  std::vector<GridFunc*> grayScottOp {&lapU, &lapV};
  DiffEq equation(pd,grayScottOp,icGrayScott,rhsGrayScott);
\end{lstlisting}
where the two separate Laplace-Beltrami operators are weighted by the different diffusivities of the $u$ and $v$ components in the Gray-Scott equations. These four lines, along with reaction terms and initial conditions, fully encode equation \eqref{eq:rdsys}.

\subsubsection{Grid functions}
\label{sec:running:user:grid}
The class \texttt{GridFunc} provides the user with a way to define their own operators without needing the specific details of the computational domain or the distributed matrix that the operator will eventually be represented by. In a two dimensional embedding space, the discrete Laplace-Beltrami operator shown in equation \eqref{eq:lapBel} can be obtained from the semidiscrete form
\begin{align}
  \Delta_{\mathcal{S}}^h u(x,y) = -\frac{4}{h^2} & u(x,y) + \nonumber \\
  \frac{1}{h^2} & \left[ u\left(cp_{\mathcal{S}}(x+h,y)\right) + u\left(cp_{\mathcal{S}}(x-h,y)\right) + \right. \\
    & \left. u\left(cp_{\mathcal{S}}(x,y+h)\right) + u\left(cp_{\mathcal{S}}(x,y-h)\right)\right]. \nonumber
  \label{eq:semiDisc}
\end{align}
Grid functions such as this are specified by three pieces of information: the stencil of nodes needed, the weight of each node in the stencil, and whether a node should be extended from the surface or not. The stencils for the extension operator are held internally, so this final item is simply a vector of booleans, and the full stencil including all nodes needed only for extension does not need to be specified.

Basic objects of the \texttt{GridFunc} class are formed from three \texttt{std::vector} objects holding these pieces of information. For ease, there are built-in grid functions for the identity, second-order accurate Laplace-Beltrami, and fourth-order accurate Laplace-Beltrami operators. Additionally, \texttt{GridFunc} objects can be added and composed to yield more complicated grid functions. For instance, the shifted Poisson equation \eqref{eq:shiftPoisson} is written in the \texttt{DDCPPoisson} example as:
\begin{lstlisting}[language=C++]
  GridFunc laplacian("LaplacianSecondOrder",
                     pd.mesh.dim,pd.mesh.delta,1.0);
  GridFunc shift("Identity",
                 pd.mesh.dim,pd.mesh.delta,1.0);
  GridFunc shiftLap = shift - laplacian;
\end{lstlisting}
where the final line encodes the total left-hand side operator.

The biharmonic operator can be written as:
\begin{lstlisting}[language=C++]
  GridFunc laplacian("LaplacianSecondOrder",
                     pd.mesh.dim,pd.mesh.delta,1.0);
  GridFunc biharmonic = laplacian*laplacian;
\end{lstlisting}
which is much simpler than writing the entire biharmonic stencil. Crucially, building these grid functions avoids doing significantly more expensive operations on large distributed matrices after their assembly.

\section{Running the software}
\label{sec:running}
The DD-CPM software acts primarily as a library. Each equation to be solved is specified in a brief driver code which calls the DD-CPM library. A number of sample driver programs are included to demonstrate how to interact with the library, and how to solve a few representative equations. A small plotting script is included in the \texttt{python} directory nested under the DD-CPM root directory. The relevant plotting command to visualize the solution from each example is given. The DD-CPM library and sample driver programs generate output both to the terminal and in the form of HDF5 data files. The interpretation and use of these outputs is discussed. With all of the sample driver programs as reference, the section is concluded with a guide towards specifying user defined equations.

\subsection{Included examples}
\label{sec:running:ex}
The DD-CPM library comes with several example programs to demonstrate usage of the library and the structure of typical driver codes. The example programs do not rely on any particular choice of surface or level of parallelism. The example programs are
\begin{itemize}
  \item \texttt{DDCPPoisson.ex:} Solves the shifted Poisson equation, demonstrating stationary solvers and Laplace-Beltrami operators.
  \item \texttt{DDCPBiharmonic.ex:} Solves the shifted biharmonic equation, demonstrating stationary solvers and composition of operators.
  \item \texttt{DDCPHeat.ex:} Solves the heat equation, demonstrating implicit time stepping.
  \item \texttt{DDCPFitzhughNagumo.ex:} Solves the Fitzhugh-Nagumo equation, demonstrating implicit/explicit time stepping and multi-component equations.
  \item \texttt{DDCPGrayScott.ex:} Solves the Gray-Scott equation.
  \item \texttt{DDCPSchnackenberg.ex:} Solves the Schnackenberg equation, and includes a user defined surface.
\end{itemize}
After compilation, the binaries for these examples can be found in the \texttt{bin} directory, and the example source codes can be found in the \texttt{examples} directory.

Running each of the example programs follows the same format:
\begin{lstlisting}[language=bash]
  mpiexec -n <nprocs> bin/DDCPExample.ex \
  -infile <path/to/inputfile> \
  <additional options>
\end{lstlisting}
where \texttt{DDCPExample} is replaced by one of the above example programs, \texttt{<nprocs>} is replaced by the number of processes you want to use for the run, and \texttt{<path/to/inputfile>} is the path to a file defining a surface and various default values. Many options can also be set on the command line to override values supplied by the input file. A full list of the available options can be obtained by running the example with the flag \texttt{-help}. The \texttt{-help} flag also gives several sample calls, specific to each example, with interesting options set to familiarize the user.

For examples using the included preconditioners (via flag \texttt{-pc\_type ddcpm}), it is required that the number of subdomains be divisible by the number of processes used. This may be changed depending on the resources available on your system. The subdomain count can easily be changed in all examples (via flag \texttt{-mesh\_nparts <n>}) keeping in mind that the number of subdomains must be divisible by the number of processes used. 

There is also a small plotting tool included in the \texttt{DD-CPM/python} directory. In the following examples, one call to the solver is shown alongside a command to generate a plot of the solution. All of the following commands are issued from inside the \texttt{DD-CPM/build} directory.

\subsubsection{DDCPPoisson}
\label{sec:ex:poisson}
The shifted Poisson equation, equation \eqref{eq:shiftPoisson} above, makes an excellent benchmark problem to compare various linear solvers and preconditioners. Similarly, it is a useful problem to showcase the flexibility available in the DD-CPM software. This equation is solved on a torus with a grid resolution of $h=1/300$ using 64 processes and a variety of linear solvers.

In all cases, the software is invoked as:
\begin{lstlisting}[language=bash]
  mpiexec -n 64 bin/DDCPPoisson.ex \
  -infile inputFiles/torus.icpm -mesh_res 300 \
  -pp_out_all -mesh_npoll 800 \
  -log_view -ksp_converged_reason \
  < additional options >
\end{lstlisting}
where \texttt{< additional options >} is set according to the desired linear solver, and consists mostly of standard PETSc flags. The flag \texttt{-mesh\_res 300} specifies the resolution, and for this surface results in a problem with \num{4090560} unknowns. The flags \texttt{-pp\_out\_all} and \texttt{-mesh\_npoll 800} control output the to polling surface, which is discussed in more detail below.

The PETSc geometric/algebraic multigrid (GAMG) \cite{petsc-web-page} solver can be specified by setting the additional options as: 
\begin{lstlisting}[language=bash]
  -mesh_nparts 1 -pc_type gamg -ksp_type richardson \
  -pc_gamg_sym_graph true -pc_gamg_agg_nsmooths 2
\end{lstlisting}
Alternatively, GAMG can be used as a preconditioner for GMRES by changing the flag \texttt{-ksp\_type richardson} to \texttt{-ksp\_type gmres}, or omitting it entirely as GMRES is the default choice within the DD-CPM software.

The MUMPs parallel sparse direct solver \cite{MUMPS01,MUMPS02} can be used by instead setting the additional options to:
\begin{lstlisting}[language=bash]
  -mesh_nparts 1 -ksp_type preonly -pc_type lu \
  -pc_factor_mat_solver_type mumps
\end{lstlisting}

The DD-CPM software includes custom DD preconditioners as discussed in Section \ref{sec:rev:dd}. The solver can be set to GMRES using these preconditioners by setting the additional options as:
\begin{lstlisting}[language=bash]
  -mesh_nparts 64 -mesh_nover 4 \
  -petscpartitioner_type parmetis \
  -pc_type ddcpm -pc_ddcpm_robfo true
\end{lstlisting}
This partitions the mesh into 64 subdomains using the ParMETIS \cite{ParMETIS} mesh partitioner, though any partitioner included in the Petsc installation can be used. These subdomains are overlapped with a width of 4, and each process is responsible for the solution of one subproblem at each iteration. The flag \texttt{-pc\_type ddcpm} enables the included DD preconditioner, and the subsequent flag enables the use of Robin transmission conditions. The parameters in the Robin transmission conditions deafult to $\alpha = 4$ and $\alpha^\times = 20$, and can be overridden with the flags \texttt{-pc\_ddcpm\_osm\_alpha <alpha>} and \texttt{-pc\_ddcpm\_osm\_alpha\_cross <alpha\_cross>}. More information regarding these transmission conditions and the parameters therein can be found in \cite{May:DD25,May:SISc}. Finally, we note that the subproblems are solved using a one level incomplete LU (iLU) factorization by default. This, and all other options regarding the solution of the subproblems, can be overridden by using the \texttt{-sub\_} prefix (e.g. \texttt{-sub\_pc\_type lu}). We emphasize that the parameter choices specified in the options above are almost certainly not optimal, and have instead been set somewhat naively to demonstrate the flexibility of the software.

The mesh resolution can be doubled from $h=1/300$ to $h=1/600$ with the flag \texttt{-mesh\_res 600}, yielding a global system with \num{16370528} unknowns. The Richardson/GAMG and GMRES/GAMG solvers can be used as is at this higher resolution. Our novel DD preconditioners can also be used, but require minor modification. For problems of this size, the use of incomplete LU (iLU) factorization can occasionally lead to stagnation of the solver. This stagnation effect arises from the interaction between iLU factorization of the local systems and the Robin transmission conditions and is not yet fully understood, although large system size has been found to be a decisive factor. Fortunately, the standard RAS preconditioner remains well behaved for large problems; Dirichlet transmission conditions are the default option, so this may be specified by omitting the flag \texttt{-pc\_ddcpm\_robfo true}. Finally, note that the provided DD preconditioners are defined only for a single level, and their efficacy starts to fall off on problems of the size considered here. The inclusion of a coarse grid correction to alleviate this effect will be investigated in future releases.



\subsubsection{DDCPBiharmonic}
This example solves the shifted biharmonic equation
\begin{equation}
  \left(c + \Delta_\mathcal{S}^2\right)u = f({\bf x}),
  \label{eq:biharm}
\end{equation}
where, as before, $c$ is a positive constant. The forcing function is given in spherical coordinates as
\begin{equation}
  f(\theta,\phi) = 5000\left(\theta^2-\frac{\pi}{2}\theta\right)\cos(7\phi),
  \label{eq:bhForce}
\end{equation}
where here, and for any uses of spherical coordinates in the remainder, $\phi\in[0,2\pi]$ is the azimuthal coordinate and $\theta\in[0,\pi]$ is the polar coordinate.

Note that the included domain decomposition preconditioners do not support biharmonic equations. Instead, one can either use preconditioners included in PETSc, or use a parallel direct solver. We take the latter approach here. Equation \eqref{eq:biharm} can be solved on a sphere with 64 processes and the MUMPS \cite{MUMPS01,MUMPS02} parallel sparse direct solver by calling:
\begin{lstlisting}[language=bash]
  mpiexec -n 64 bin/DDCPBiharmonic.ex \
  -infile inputFiles/sphere.icpm \
  -mesh_res 200 -mesh_nparts 1 -ksp_type preonly \
  -pc_type lu -pc_factor_mat_solver_type mumps
\end{lstlisting}
where the flags \texttt{-mesh\_res 200} and \texttt{-mesh\_nparts 1} set a grid spacing of $1/200$ and inform the library that no mesh partitioning is necessary. The remaining flags are all standard PETSc options to enable the MUMPS solver.

The grid function defining this operator can be obtained directly from the operator for the Laplacian by composition. Crucially, this is done before the global operator is constructed to avoid expensive matrix-matrix multiplication. This composition is quite simple, and is discussed in more detail in Section \ref{sec:running:user:grid}.

\begin{figure}
  \centering
  \includegraphics[width=0.8\linewidth]{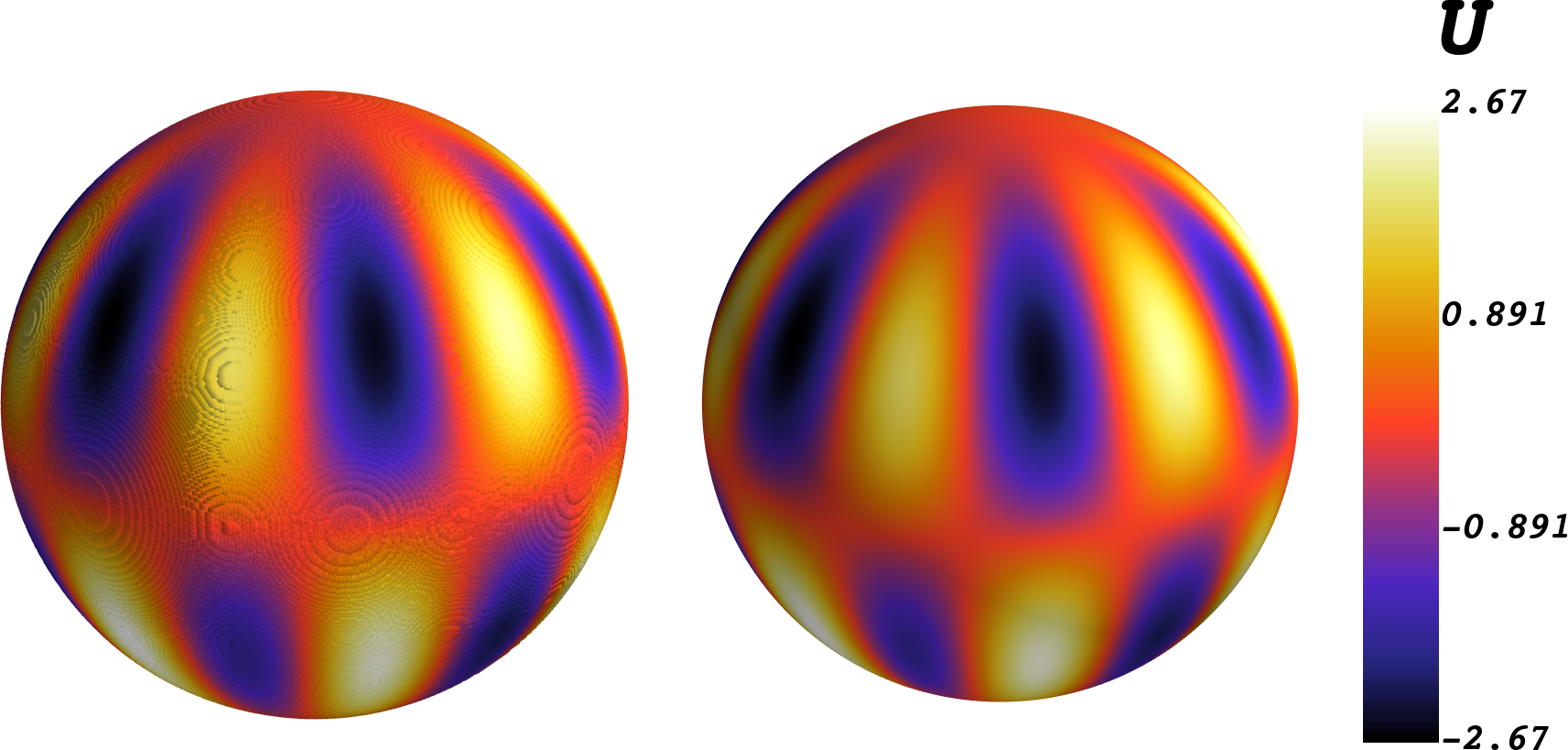}
  \caption{The left and right panels show the solution of the shifted biharmonic equation \eqref{eq:biharm} on the point cloud and polling surface respectively. The grid spacing is $h=1/200$, and the forcing function is given in equation \eqref{eq:bhForce}.}
  \label{fig:biharmonic}
\end{figure}

The solution can be visualized on the CPM point cloud by calling the included plotting tool as:
\begin{lstlisting}[language=bash]
  python ../python/cpmplot.py \
  data/Sphere200_2_1_3/Sphere000000.h5 cloud sol
\end{lstlisting}
or over the surface itself by:
\begin{lstlisting}[language=bash]
  python ../python/cpmplot.py \
  data/Sphere200_2_1_3/Sphere000000.h5 poll sol
\end{lstlisting}
where \texttt{poll} indicates that the solution should be displayed on the \textit{polling surface} discussed in more detail below. The solution can be seen in Figure \ref{fig:biharmonic}.

\subsubsection{DDCPHeat}
This example solves the heat equation
\begin{equation}
  \frac{\partial u}{\partial t} - \Delta_\mathcal{S}u = f({\bf x},t),
  \label{eq:heat}
\end{equation}
where the time-dependent forcing function is given in spherical coordinates as
\begin{equation}
  f(\theta,\phi,t) = \cos(t)\sin(3\phi).
  \label{eq:htForce}
\end{equation}
The initial condition is $u({\bf x},0)=0$.

This can be solved on a sphere using 4 processes and our ORAS preconditioner over 12 subdomains by calling:
\begin{lstlisting}[language=bash]
  mpiexec -n 4 bin/DDCPHeat.ex \
  -infile inputFiles/sphere.icpm \
  -mesh_res 30 -mesh_nparts 12 -mesh_nover 4 \
  -petscpartitioner_type parmetis \
  -pc_type ddcpm -pc_ddcpm_robfo true \
  -time_final 0.2 -ts_monitor
\end{lstlisting}
The default time integrator is the PETSc TSARKIMEX method, and the linear systems that arise are solved with GMRES supplemented by a 12 subdomain DD-CPM ORAS preconditioner. As above, the flags \texttt{-petscpartitioner\_type parmetis -pc\_type ddcpm -pc\_ddcpm\_robfo true} indicate that the mesh should be partitioned using ParMETIS, and that Robin transmission conditions should be used with their default parameters of $\alpha=4$ and $\alpha^{\times}=20$. Again, we note that iLU factorization is used in the solution of all subproblems, which is generally useful for transient problems.

\begin{figure}
  \centering
  \includegraphics[width=0.8\linewidth]{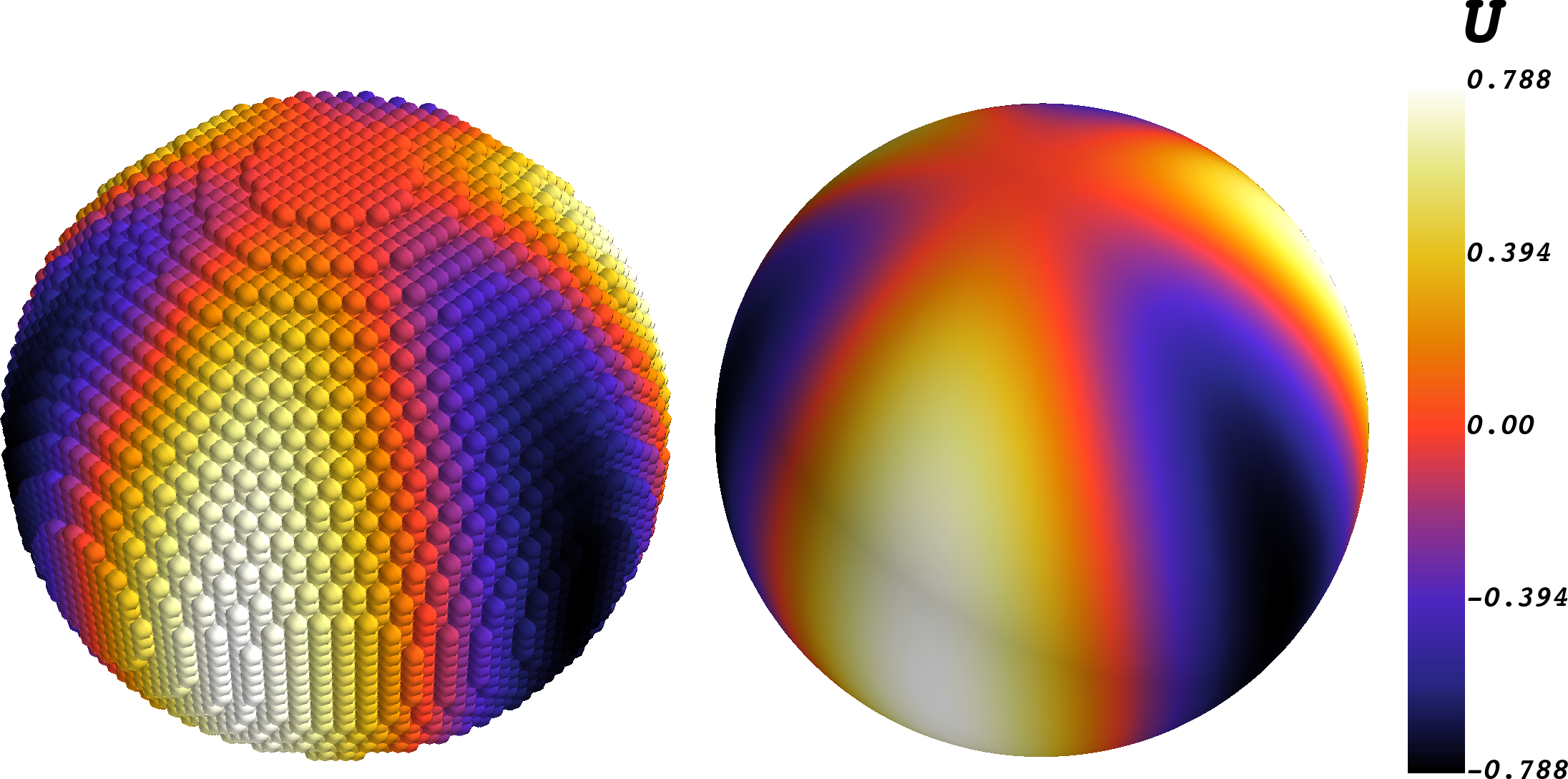}
  \caption{The left and right panels show the solution of the heat equation \eqref{eq:heat} on the point cloud and polling surface respectively. The grid spacing is $h=1/30$, the final time is $t=0.2$, and the forcing function is given in equation \eqref{eq:htForce}.}
  \label{fig:heat}
\end{figure}

The solution can be visualized on the CPM point cloud by calling the included plotting tool as:
\begin{lstlisting}[language=bash]
  python ../python/cpmplot.py \
  data/Sphere30_2_12_4/Sphere000006.h5 cloud sol
\end{lstlisting}
or over the polling surface by:
\begin{lstlisting}[language=bash]
  python ../python/cpmplot.py \
  data/Sphere30_2_12_4/Sphere000006.h5 poll sol
\end{lstlisting}
both of which can be seen in Figure \ref{fig:heat}.

\subsubsection{DDCPFitzhughNagumo}
This example solves the Fitzhugh-Nagumo equations
\begin{equation}
  \begin{cases}
    \frac{\partial u}{\partial t} - D_u\Delta_\mathcal{S}u = (0.1-u)(u-1)u - v, \\
    \frac{\partial v}{\partial t} - D_v\Delta_\mathcal{S}v = 0.01(0.5u-v)
  \end{cases}
  \label{eq:fitznag}
\end{equation}
where the diffusivities of each species are set as $D_u=10^{-4}$ and $D_u=10^{-7}$. The initial condition is chosen to be
\begin{align}
  u({\bf x},0) &= \begin{cases} 1,~& x>0,y>0,z>0 \\ 0,~&{\rm otherwise} \end{cases} \\
  v({\bf x},0) &= \begin{cases} 1,~& x>0,y<0,z>0 \\ 0,~&{\rm otherwise} \end{cases}.
  \label{eq:fnICs}
\end{align}

Time dependent problems require more setup, though most of the flags are self explanatory. This can be solved on a triangulated surface using 12 processes and the PETSc built-in block Jacobi preconditioner:
\begin{lstlisting}[language=bash]
  mpirun -n 12 bin/DDCPFitzhughNagumo.ex \
  -infile inputFiles/triang.icpm \
  -mesh_res 40 -pc_type bjacobi \
  -time_final 250 -ts_monitor
\end{lstlisting}
The triangulated surface is supplied in the form of a PLY file which is specified in the \texttt{triang.icpm} input file. The default option is the Stanford Bunny \cite{bunny}, though several other PLY files are included in the \texttt{DD-CPM/plyFiles} directory.

\begin{figure}
  \centering
  \includegraphics[width=0.8\linewidth]{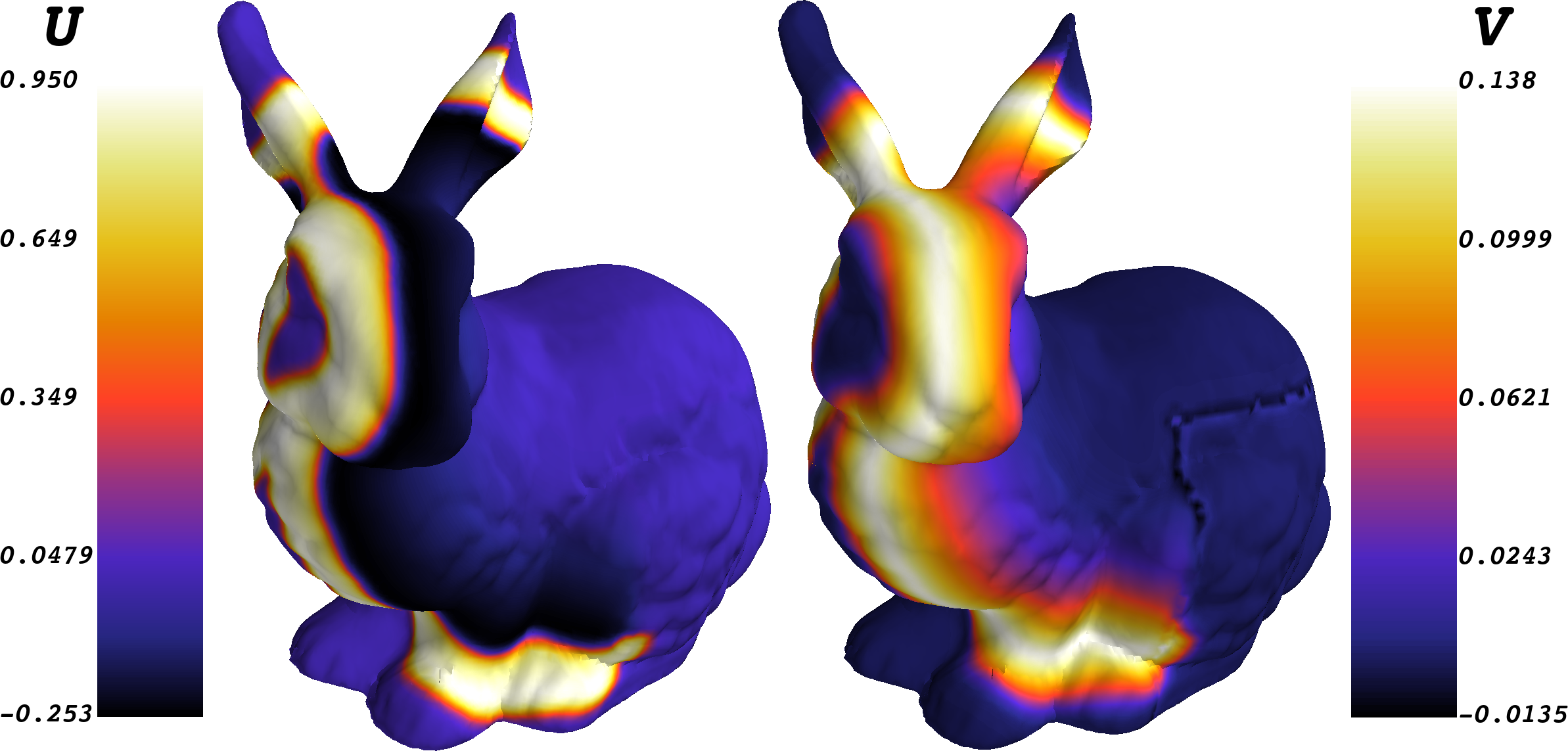}
  \caption{The left and right panels show the $u$ and $v$ components of the solution to the Fitzhugh-Nagumo equation \eqref{eq:fitznag} at time $t=250$. The initial conditions are set according to equation \eqref{eq:fnICs}, and the grid spacing used is $h=1/40$.}
  \label{fig:fitznag}
\end{figure}

The $u$ component of the solution can be visualized by calling:
\begin{lstlisting}[language=bash]
  python ../python/cpmplot.py \
  data/Triangulated40_2_12_3/Triangulated000176.h5 \
  poll sol 0
\end{lstlisting}
and similarly the $v$ component of the solution can be visualized by calling:
\begin{lstlisting}[language=bash]
  python ../python/cpmplot.py \
  data/Triangulated40_2_12_3/Triangulated000176.h5 \
  poll sol 1
\end{lstlisting}
each of which can be seen in Figure \ref{fig:fitznag}. Note that the solution here is not shown on the point cloud, but rather on the surface itself. This is specified by the \texttt{poll} option. When the DD-CPM library produces output files it will try to project the solution back onto the surface for later visualization. This requires that some triangulation of the surface is known. Of course, this triangulation is \textit{not} used to solve the equations, but is needed to produce these visualizations. The solution on the point cloud can be seen by replacing \texttt{poll} with \texttt{cloud} in the above plotting calls.

\subsubsection{DDCPGrayScott}
\label{sec:ex:grayscott}
This example solves the Gray-Scott equations
\begin{equation}
  \begin{cases}
    \frac{\partial u}{\partial t} - D_u\Delta_\mathcal{S}u = (1-u)F - uv^2\\
    \frac{\partial v}{\partial t} - D_v\Delta_\mathcal{S}v = uv^2 - (F+K)v
  \end{cases}
  \label{eq:grayscott}
\end{equation}
where the diffusivities of each species are chosen as $D_u = 6\times 10^{-5}$ and $D_u = 3\times 10^{-5}$, the feed rate is set at $F=0.03$, and the kill rate is given in spherical coordinates by the function
\begin{equation}
  k(r,\theta,\phi) = 0.06 + 0.001\cos(\phi)
  \label{eq:gsPars}
\end{equation}
where $\phi$ is the azimuthal angle. The initial condition is chosen to be
\begin{align}
  u(r,\theta,\phi,0) &= 1-e^{-5\sin^2(5(\phi+0.5))} \\
  v(r,\theta,\phi,0) &= e^{-5\sin^2(5(\phi-0.5))}.
  \label{eq:gsICs}
\end{align}

As in the previous example, the nonlinear coupling terms are advanced through time explicitly. To solve the Gray-Scott equation on a Dupin cyclide with 128 processes one could call:
\begin{lstlisting}[language=bash]
  mpiexec -n 128 bin/DDCPGrayScott.ex \
  -infile inputFiles/triang.icpm -mesh_res 300 \
  -mesh_nparts 128 -mesh_nover 4 \
  -petscpartitioner_type parmetis \
  -pc_type ddcpm -pc_ddcpm_robfo true \
  -time_final 15000 -pp_plotfreq 100 -ts_monitor
\end{lstlisting}
This generates a mesh with \num{3739898} active nodes globally. The diffusion terms are treated implicitly, and the linear systems that arise are solved with GMRES preconditioned by a 128 subdomain DD-CPM ORAS preconditioner. The time stepper defaults to using the PETSc TSARKIMEX3 method, and the present flags set a final time of $T=15,000$. Finally, the flag \texttt{-pp\_plotfreq 100} informs the post-processor to write a data file every 100 time steps.

\begin{figure}[hbtp]
  \centering
  \includegraphics[width=0.95\linewidth]{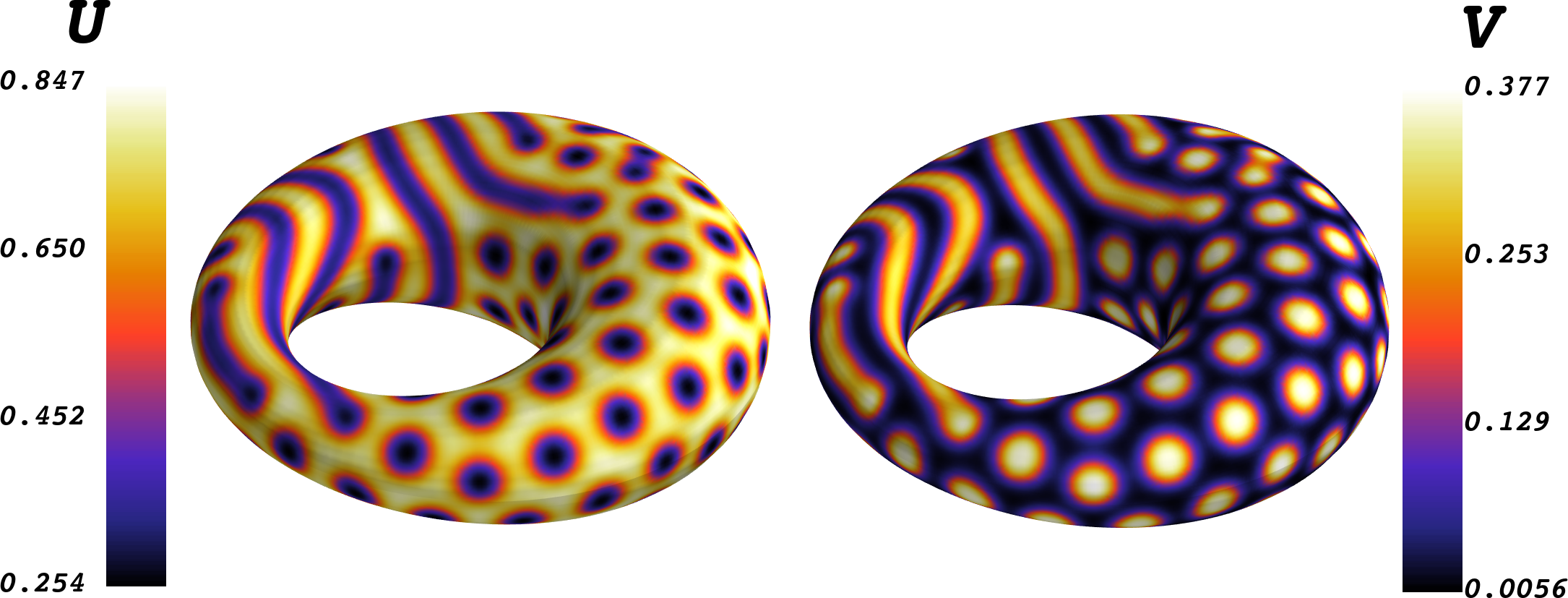}
  \caption{The left and right panels show the $u$ and $v$ components of the solution to the Gray-Scott equation \eqref{eq:grayscott} at time $t=15000$. The initial conditions and parameters are set according to equations \eqref{eq:gsICs} and \eqref{eq:gsPars}, and the grid spacing used is $h=1/400$.}
  \label{fig:grayscott}
\end{figure}

These solutions can be visualized by calling the included plotting tool.    The $u$ component of the solution can be visualized by calling
\begin{lstlisting}[language=bash]
  python ../python/cpmplot.py \
  data/Triangulated300_2_128_4/Triangulated000510.h5 \
  poll sol 0
\end{lstlisting}
while  the $v$ component can be visualized via
\begin{lstlisting}[language=bash]
  python ../python/cpmplot.py \
  data/Triangulated300_2_128_4/Triangulated000510.h5 \
  poll sol 1
\end{lstlisting}
In each call \texttt{poll} specifies that the polling surface should be used for visualization. The solutions can be seen in Figure \ref{fig:grayscott}. To see the solution on the CPM point cloud one would need to replace \texttt{poll} by \texttt{cloud} in the visualization command.

\subsubsection{DDCPSchnackenberg} The final example solves the Schnackenberg equations
\begin{equation}
  \begin{cases}
    \frac{\partial u}{\partial t} - D_u\Delta_\mathcal{S}u = \delta_1 - u + u^2v + \lambda v, \\
    \frac{\partial v}{\partial t} - D_v\Delta_\mathcal{S}v = \delta_2 - u^2v - \lambda v
  \end{cases}
  \label{eq:schnack}
\end{equation}
where the parameters are chosen to be $D_u = 0.001$, $D_v = 0.1$, $\delta_1 = 0.003$, $\delta_2 = 0.7$, and $\lambda = 0.06$ inspired in part by \cite{AlNoufaey:Schnack}. The initial conditions are set as
\begin{align}
  u &= 0.02 + \mathcal{U}\sin(5x)\cos(3y), \\
  v &= 4 + \mathcal{V}\sin(5xy)\sin(5x^3y)
  \label{eq:schnackIcV}
\end{align}
where $\mathcal{U}\sim\mathcal{N}(0,0.01)$ and $\mathcal{V}\sim\mathcal{N}(0,1)$ are normally distributed random values.

This example also demonstrates how to implement user defined surfaces. Inside the driver code \texttt{DDCPSchnackenberg.cpp} there are functions defining the closest point function, a distance function, and a surface normal function, which are all handed to the library. In this case these functions correspond to a unit hemisphere.

The boundary conditions along the circular boundary are mixed between homogeneous Neumann and homogeneous Dirichlet. Taking $\phi$ as the azimuthal angle, any point that maps to the boundary where $\sin(3\phi) < -0.5$ is mirrored to enforce second order accurate homogeneous Dirichlet boundary conditions. The remainder of the circular boundary utilizes the natural homogeneous Neumann boundary conditions. See \cite{CBM:Eig} for details of the Dirichlet boundary formulation.

\begin{figure}
  \centering
  \includegraphics[width=0.8\linewidth]{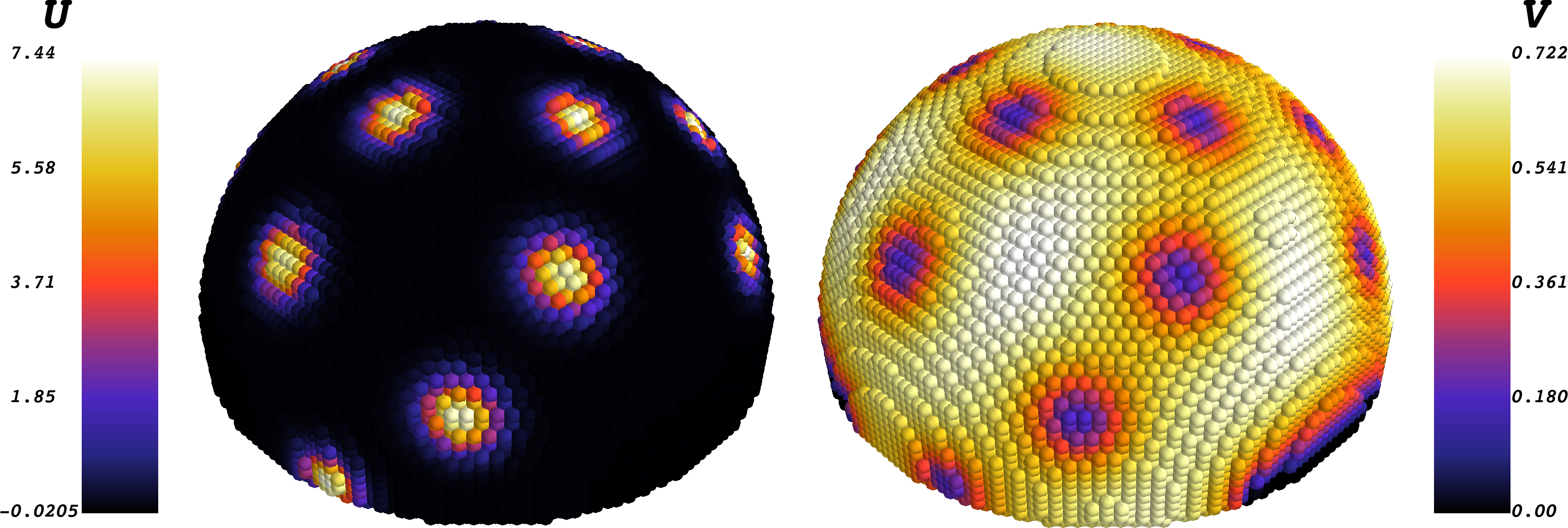}
  \caption{The left and right panels show the $u$ and $v$ components of the solution to the Schnackenberg equation \eqref{eq:schnack} at time $t=50$. The initial conditions are given by equation \eqref{eq:schnackIcV}, the grid spacing is $h=1/30$, and the parameters used are $D_u = 0.001$, $D_v = 0.1$, $\delta_1 = 0.003$, $\delta_2 = 0.7$, and $\lambda = 0.06$.}
  \label{fig:schnack}
\end{figure}

This equation can be solved on the user defined surface using 12 processes and our ORAS preconditioner over 12 subdomains by calling:
\begin{lstlisting}[language=bash]
  mpirun -n 12 bin/DDCPSchnackenberg.ex \
  -infile inputFiles/user.icpm -cp_pc_ras -mesh_res 30 \
  -mesh_nparts 12 -mesh_nover 4 \
  -petscpartitioner_type parmetis \
  -pc_type ddcpm -pc_ddcpm_robfo true \
  -time_final 100 -ts_monitor
\end{lstlisting}

The $u$ component of the solution can be visualized by calling:
\begin{lstlisting}[language=bash]
  python ../python/cpmplot.py \
  data/User30_2_12_4/User000832.h5 cloud sol 0
\end{lstlisting}
and similarly the $v$ component of the solution can be visualized by calling:
\begin{lstlisting}[language=bash]
  python ../python/cpmplot.py \
  data/User30_2_12_4/User000832.h5 cloud sol 1
\end{lstlisting}
where each can be seen in Figure \ref{fig:schnack}. Note that each of these calls visualize the solution over the point cloud. Currently there is no support for polling surfaces corresponding to user defined surfaces.

\section{Profiling}
\label{sec:num}
The DD-CPM library yields great parallel scalability for moderate and large problems using hundreds of cores spread across a few nodes of a compute cluster. To demonstrate this scalability we focus on the shifted Poisson equation due to its ubiquity and relative simplicity. We remark before presenting the results that mesh partitioning can be done by any of the PETSc supported graph partitioners, but does have a minor serial component afterward. This does present an early bottleneck when running large problems. However, this only occurs once in any run and does not pose any issue for problems with millions of active grid points. This bottleneck will be removed in a future release of the library. All tests were performed on whole nodes of the Graham compute cluster managed by Compute Canada. Each node has two 16-core Intel E5-2683 v4 Broadwell 2.1 GHz processors and 125G of memory.

Each run is broken into three major phases: mesh construction and partitioning, global and local operator construction, and the GMRES solve. As discussed above, the meshing phase is not expected to scale due to the serial bottleneck there. The time required for each phase is recorded using the logging features of PETSc. Each test of scalability is performed using $N_P=64,~128,~256$ processes, corresponding to 2, 4, and 8 full compute nodes on Graham respectively.

To assess the strong scalability, where a problem of fixed size is solved with progressively greater amounts of parallelism, we consider the torus described in Section \ref{sec:ex:poisson} with a grid spacing of $h=\frac{1}{300}$ and tri-quadratic interpolation. This combination yields a problem with \num{4090560} unknowns. The linear system is solved with GMRES using the DD-CPM RAS preconditioner with an overlap width of $N_O=4$. All local problems inside the RAS preconditioner are solved inexactly using a one level incomplete LU factorization.

\begin{table}[hbt]
  \centering
  \begin{tabular}{@{}cccc@{}}
    \toprule
    Phase                 & $N_P=64,~N_S=64$ & $N_P=128,~N_S=128$ & $N_P=256,~N_S=256$ \\ \midrule
    Meshing               & 33               & 33              & 34                 \\
    Operators             & 20               & 11              & 8                  \\
    GMRES                 & 33               & 22              & 12                 \\
    Total                 & 85               & 65              & 55                 \\ \midrule
    Iterations            & 391              & 496             & 481                \\ 
    Time per iteration    & 0.084            & 0.044           & 0.025              \\ \bottomrule
  \end{tabular}
  \caption{The time in seconds is shown for each phase in the solution of the shifted Poisson equation on a torus with $h=300^{-1}$ and \num{4090560} unknowns. The required number of iterations and time per iteration is listed in the lower portion of the table. Here $N_S=N_P$ is fixed and good strong scalability is observed apart from the meshing phase.}
  \label{tab:sclStrEq}
\end{table}

The number of subdomains used must be divisible by the number of processes, giving two interesting options. First, Table \ref{tab:sclStrEq} presents results from setting $N_S=N_P$. For $N_S=N_P=64$ each disjoint subdomain contains approximately \num{64000} active nodes, and each overlapping subdomain contains approximately \num{73000} nodes. For $N_S=N_P=128$ each disjoint subdomain contains approximately \num{32000} active nodes, and each overlapping subdomain contains approximately \num{38000} nodes. Finally, for $N_S=N_P=256$ the disjoint and overlapping subdomains contain approximately \num{16000} and \num{20000} active nodes respectively. In this case the preconditioner is stronger for smaller process counts, but more expensive to apply due to the larger local problems. The operator construction and GMRES solve scale quite well, though the disparity in the size of the local problems is visible when comparing $N_P=64$ to $N_P=128$.

\begin{table}[hbt]
  \centering
  \begin{tabular}{@{}cccc@{}}
    \toprule
    Phase                 & $N_P=64,~N_S=256$ & $N_P=128,~N_S=256$ & $N_P=256,~N_S=256$ \\ \midrule
    Meshing               & 33                & 34               & 34                 \\
    Operators             & 21                & 12               & 8                  \\
    GMRES                 & 46                & 22               & 12                 \\
    Total                 & 99                & 67               & 55                 \\ \midrule
    Iterations            & 554               & 480              & 481                \\
    Time per iteration    & 0.083             & 0.046            & 0.025              \\ \bottomrule
  \end{tabular}
  \caption{The time in seconds is shown for each phase in the solution of the shifted Poisson equation on a torus with $h=300^{-1}$ and \num{4090560} unknowns. Here $N_S=256$ is fixed independent of the number of processes (compare to Table \ref{tab:sclStrEq}). Excellent strong scalability is observed apart from the meshing phase.}
  \label{tab:sclStrFix}
\end{table}

Alternatively, the number of subdomains can be fixed at $N_S=256$ irrespective of the number of processes. This ensures that all of the local problems are roughly the same size, and that the preconditioner being applied is the same in all cases. In this case, each disjoint subdomain contains approximately \num{16000} active nodes, and each overlapping subdomain contains approximately \num{20000} nodes. Note that for $N_P=64$ and $N_P=128$ each process is responsible for constructing and factoring multiple local operators. 
As shown in Table \ref{tab:sclStrFix}, the now weakened preconditioner (for the smaller process counts) leads to longer GMRES solution times. Excellent strong scalability is observed apart from the meshing phase.

\begin{table}[hbt]
  \centering
  \begin{tabular}{@{}cccc@{}}
    \toprule
    Phase                 & $N_P=64,~N_S=64$ & $N_P=128,~N_S=128$ & $N_P=256,~N_S=256$ \\ \midrule
    Meshing               & 33               & 65              & 139                \\
    Operators             & 20               & 21              & 25                 \\
    GMRES                 & 33               & 77              & 121                \\
    Total                 & 85               & 164             & 288                \\ \midrule
    Iterations            & 391              & 877             & 1438               \\
    Time per iteration    & 0.084            & 0.088           & 0.084              \\ \bottomrule
  \end{tabular}
  \caption{The time in seconds is shown for each phase in the solution of the shifted Poisson equation on a torus on a sequence of finer meshes. Resolutions of $h=300^{-1}$, $h=420^{-1}$, and $h=595^{-1}$, were used to keep the number of unknowns per process roughly constant.  Overall solution times increase somewhat due to increasing times in  the meshing phase (which contains a serial component) and an increase in the number of GMRES iterations (which is as expected due to increasing subdomain size).  Excellent weak scalability is observed in the time per iteration.  
 }
  \label{tab:sclWkEq}
\end{table}

Weak scalability is assessed by fixing $N_S=N_P$, and varying the resolution of the problem such that the sizes of the disjoint partitions are roughly constant. For $N_S=N_P=64$ the resolution is $h=\frac{1}{300}$ for a total of \num{4090560} unknowns with disjoint and overlapping subdomains of approximate size \num{64000} and \num{73000} respectively. For $N_S=N_P=128$ the resolution is increased to $h=\frac{1}{420}$ for a total of \num{8060600} unknowns, and for $N_S=N_P=256$ the resolution is further increased to $h=\frac{1}{595}$ for a total of \num{16082068} unknowns. As can be seen in Table \ref{tab:sclWkEq}, the operator construction phase scales almost perfectly. The time required for the GMRES solution slowly climbs, which is to be expected since the included preconditioner incorporates only a single level. This is demonstrated concisely by the growing number of required iterations to convergence.   Excellent weak scalability is observed in the time per iteration.

\section{Conclusion}
\label{sec:conc}
The DD-CPM library provides a general purpose implementation of the closest point method with an emphasis on distributed memory parallelism. The DD-CPM library is particularly useful for surface intrinsic reaction-diffusion equations, where the included DD preconditioners can efficiently handle the stiff diffusion terms implicitly. The DD-CPM library also handles domains in flat space with otherwise difficult boundaries, and with some development could be applied to coupled bulk-surface equations. The software could be similarly extended to handle fully implicit discretizations in time, with the given preconditioners aiding in the solution of the linear systems that arise inside the Newton iterations. Additionally, the DD-CPM library could make an effective back end for eigenvalue problems due to the interoperability of PETSc and SLEPc.


\section*{Declarations}
\subsection*{Funding}
  The authors gratefully acknowledge the financial support of NSERC Canada (RGPIN 2016-04361 and RGPIN 2018-04881). This research was enabled in part by support provided by ACEnet (ace-net.ca) and Compute Canada (www.computecanada.ca).
\subsection*{Data Sharing}
All data generated as part of this study is available upon request or can be generated using the described library and the calling sequences listed in the paper. 

\subsection*{Conflict of Interest}
The authors declare that they have no conflict of interest.

\bibliography{references}

\end{document}